\documentclass[12pt]{article}
\usepackage{rrfam}
\usepackage[noblocks,auth-sc]{authblk}
 3
 2

\def\bigspace{\advance\lineskip by 3pt
      \advance\baselineskip by 3pt
      \advance\lineskiplimit by 3pt}
\global\def\sectitle#1\par{\bigbreak
  \leftline{\bf #1}
  \nobreak\medskip\vskip-\parskip
  \message{#1}
  \noindent}
\global\def\ssectitle#1\par{\bigbreak\medskip
  \leftline{\typc #1}
  \nobreak\bigskip\vskip-\parskip
  \message{#1}
  \noindent}
\global\def\sssectitle#1\par{\bigbreak\bigskip
  \leftline{\typd #1}
  \nobreak\bigskip\vskip-\parskip
  \message{#1}
  \noindent}
\global\def\sssectitletwo#1#2#3\par{\bigbreak\bigskip
  \vbox{
  \leftline{\typd #1#2}
  \nobreak\vskip2truemm
  \message{#1#2}
  \leftline{\typd \phantom{#1}#3}
  \message{#3}
  }
  \nobreak\bigskip\vskip-\parskip
  \noindent}

\newcounter{saveeqn}

\def\mod{\mathop{\rm mod}}

\def\dim{\mathop{\rm dim}}
\def\Fix#1{\mathop{\rm Fix}(#1)}

\def\det{\mathop{\rm det}}

\def\pt{\partial}

\def\inv{^{-1}}

\def\Meng#1#2{\left\{#1\;\Big|\;#2\right\}}
\def\intl#1#2#3{\int\limits_{#1}^{#2}{#3}}

\def\qed{\vbox{\hrule
  \hbox{\vrule\hbox to 5pt{\vbox to 8pt{\vfil}\hfil}\vrule}\hrule}}
\def\endproof{\unskip \nobreak \hskip0pt plus 1fill \qquad \qed \par}

\let\imp=\Rightarrow
\let\leftv=|
\let\rightv=|

\def \nid{\noindent }

\def\eps{\varepsilon }

\def\transversal{\kern0.5em\cap\kern-1em\stackrel{\textstyle\top}{}\kern0.25em}
\def\OO#1{\mathop{{\bf O}(#1)}}

\def\SO#1{\mathop{{\bf SO}(#1)}}

\def\DD#1{\mathop{{\D}_{#1}}}

\def\ZZ#1{\mathop{{\Z}_{#1}}}

\def\1{{\rm 1\hskip-0.9truemm l}}

\def\ifundefined#1{\expandafter\ifx\csname
  #1\endcsname\relax}
\ifundefined{textheight}
\def\newline{\par\noindent}
\else
\fi

\newtheorem{xggg}{ggg}[section]
\newtheorem{lem}[xggg]{Lemma}

\newtheorem{thm}[xggg]{Theorem}

\def\pf{{\bf Proof.} }

\author{Reiner Lauterbach}\affil{Fachbereich Mathematik\\Universit\"at
  Hamburg} \author{Paul Matthews}\affil{School of Mathematical Sciences\\University of Nottingham}
\title{Do Absolutely Irreducible Group Actions Have Odd Dimensional Fixed Point Spaces?}
\begin{document}
\maketitle
\begin{abstract}
  In his volume  \cite{Fd6} on ``Symmetry Breaking for Compact Lie
  Groups'' Mike Field quotes a private communication by Jorge Ize
  claiming that any 
  bifurcation problem with absolutely irreducible group action would lead to
  bifurcation of steady states. The proof should come from the fact that any
  absolutely irreducible representation possesses an odd dimensional fixed
  point space.\\
In this paper we show that there are many examples of groups which have
absolutely irreducible representations but no odd dimensional fixed
point space. This observation may be relevant also for some degree theoretic
considerations concerning equivariant bifurcation. Moreover we show that our
examples give rise to some interesting Hamiltonian dynamics and we show that
despite some complications we can go a long way towards doing explicit
computations and providing complete proofs. For some of the invariant theory
needed we will depend on some computer aided computations. The work presented
  here greatly benefited from the computer algebra program GAP \cite{GAP},
  which is an 
  indispensable aid for doing the required group theory computations.
\end{abstract}
\section{Introduction}
Symmetric systems may have solutions with less symmetry than the original
problem. In bifurcation theory one can have the situation that
branching of fully 
symmetric solutions leads to less symmetric solutions. Here the notion
of a solution
refers to steady states of some dynamical problem. Questions on
symmetries of solutions
come up in pattern formation and other applied problems. The general problem
can be stated as follows: let $X$ be a state space with a group action, let us
say by a compact Lie group $G$, let $P$ be a parameter space and suppose that $F:X \times P\to X$
is the right hand side of a differential equation which is equivariant with respect to this action, i.e.
\[F(gx,p)=gF(x,p)\]
for all $g\in G$, $p\in P$. If $x_0$ for given $p_0$  is a steady
state solution which is fully symmetric, 
then we have
\[F(gx_0,p_0)=0\]
for all $g\in G$. Symmetry breaking occurs if for $p$ near
$p_0$ we find steady state solutions $x(p)$ with $gx(p)\not=x(p)$ for at least
one $g\in G$. Steady state bifurcation with symmetry has a long
history, see for example
Vanderbauwhede \cite{V1}, Sattinger, \cite{Sat1}, Golubitsky et al.
\cite{GS1,GSS}, Chossat and Lauterbach \cite{ChL}. One of the main
results is the so called Equivariant Branching Lemma. It addresses the
situation when $X$ is a real Banach space and the branching comes from
a change of stability of the steady state solution where at
criticality the linearisation has a kernel which is an absolutely
irreducible representation of $G$. \\ 
Let us briefly recall these notions. If $X$ is a Banach space,
$F:X\times P\to X$ is sufficiently smooth and equivariant. Let
$F(x_0,p_0)=0$ and assume that $K=D_xF(x_0,p_0)$ is the kernel of the
linearisation of $F$ at this point. If $D_xF(x_0,p_0)$ is an
isomorphism of $X$, then locally near
$p_0$ the solution manifold can be parameterized over $p$. If we assume
the parameter space to be one-dimensional, then the solution manifold
is locally a one dimensional manifold. If we assume, that at $p_0$ the
operator $D_xF(x_0,p_0)$ is not an isomorphism, and moreover if we
assume, that $\dim\ker D_xF(x_0,p_0)>0$, then the kernel
$K=D_xF(x_0,p_0)$ is invariant under 
the group action. It is a generic property that $K$ is an
absolutely irreducible representation, i.e. that a linear map
commuting with the group action is a multiple of the identity. In such
a situation the eigenvalue $0$ can be prolonged to neighboring
parameter values and the group action will remain the same \cite{ChL}. Therefore
the multiplicity of the critical eigenvalue will not change and we
have to look at solutions in the kernel $K$. The geometry of the group
action helps to overcome some of the problems of the higher multiplicity
of the eigenvalue.    

\begin{thm}[Sattinger, Vanderbauwhede, Cicogna] \label{thm:ebl}\cite{GSS,IG}
  If $H<G$ is an isotropy subgroup with $\dim \Fix{H}=1$, then the
  Hopf   condition 
\[\sigma'(p)\not=0\]
implies bifurcation of a branch of steady states with isotropy $H$. 
\end{thm}
Here we assume that the space of parameters is one-dimensional and
$\sigma(p)$ is the curve of eigenvalues prolonging the critical
eigenvalue.
It is an open question whether the loss of stability through an
absolutely irreducible kernel always leads to a bifurcation of a
branch of steady states. In the abstract we find a strategy to prove
such a result. It relies on a slight generalization of the Equivariant
Branching Lemma.
\begin{thm} \label{thm:gen_ebl}\cite{ChL} Given a compact Lie group $G$
an $G$-equivariant bifurcation
  problem 
\[F(x,\lambda)=0\]
with an absolutely irreducible $G$-action on the kernel 
\[K=\ker D_xF(0,0)\]
of the
  linearisation. 
  If $H<G$ is an isotropy subgroup for which $\dim \Fix{H}$ is odd, then the
  Hopf   condition 
\[\sigma'(p)\not=0\]
implies bifurcation of a branch of steady states with isotropy at
least $H$. 
\end{thm}
Here as before we write $\sigma$ for the prolongation of the critical
eigenvalue. 

Then, if we can show that each group with an absolutely
irreducible group action has an isotropy subgroup with an
odd dimensional kernel, then a general result on bifurcation of equilibria
in the presence of
absolutely irreducible group actions would follow from the above
statement. As far as we are aware, this property holds for all group
actions previously considered in the context of equivariant
bifurcation theory. In particular, it is true for all absolutely
irreducible group actions on $\R^2$ and $\R^3$.
In dimension $2$ the groups acting
absolutely irreducible are $\DD{m}$, $m\ge 3$ and $\OO{2}$. These groups
all contain reflections which have a one-dimensional fixed point space.
In $\R^3$, the relevant groups are $\OO{3}$ and some of its subgroups,
and all of these groups contain a rotation with a one-dimensional
fixed point space.

In this paper,
we will show that this strategy cannot be successful, by
providing infinite series of finite groups each of which acts
absolutely irreducibly on $\R^4$, for which the only non-trivial isotropy
subgroups have two dimensional fixed point subspaces.
In section~\ref{sec:quat}, the powerful and compact quaternion
notation for actions on $\R^4$ is introduced. The subgroups of
interest and our main results are given in section~\ref{sec:series}.
The equivariant vector fields and their Hamiltonian structure 
are discussed in sections~\ref{sec:flow} and \ref{sec:Hamilton}, followed
by the proofs of the main results in section~\ref{proofs}. 
Some computational results are given in  section~\ref{sec:gap}.

\section{The group $\SO{4}$, its subgroups and quaternion notation}
\label{sec:quat}

A classification of subgroups of $\SO{4}$ and $\OO{4}$ goes back to Goursat
\cite{Gours}; some data relevant for bifurcation theory has been given by Becker and
Kr\"amer \cite{BeKr-93}. Here we use the classification of subgroups as it is
presented by Conway and Smith \cite{CS}. In this recent (and very nice) book the
quaternions are used to give a geometric way to describe the subgroups of
$\SO{3}, \OO{3}, \SO{4}$ and some others.
This quaternion notation provides a much more  compact and elegant
description of  $\SO{4}$ than the use of $4\times 4$ matrices.
Although we will use the form described by \cite{CS}, it is worth
noting that an equivalent notation was described by Felix Klein
\cite{K08,K08e}, who in turn attributes the key result to Cayley \cite{Cay1855}.

We will denote the set of unit
quaternions by $Q$. The set of pairs of such quaternions forms a
six-dimensional group, called
the spinor group and denoted by $\mbox{Spin}_4$. We get a map
\begin{equation}\label{quaternionmap}
\mbox{Spin}_4\to \SO{4}: (l,r)\mapsto [l,r]=\{x\mapsto \bar l xr\}
\end{equation}
where a vector in $\R^4$ is identified with a quaternion via
\[x=\left(
  \begin{array}{c}
x_1\\ x_2\\ x_3\\ x_4
  \end{array}\right)\quad \Leftrightarrow\quad x_1+i\, x_2+j\, x_3+k\,
x_4.\]
It is obvious that
\begin{equation}\label{series.identify}[-l,-r]=[l,r]\end{equation}
and \cite{CS} show that this is the only way that injectivity fails,
so the map is two-to-one.

In a similar way we can obtain a map $\mbox{Spin}_4\to \OO{4}\setminus\SO{4}$ by
 \[\mbox{Spin}_4\to \OO{4}\setminus \SO{4}: (l,r)\mapsto \star [l,r]=\{x\mapsto
 \bar l\bar x r\},\]
but this form will not be used in this paper.
Using the map (\ref{quaternionmap}), composition of elements of
$\SO{4}$ can be written in a natural way,
\[
[l_1, r_1] [l_2, r_2] = [l_1 l_2, r_1 r_2],
\]
since $\overline{l_1 l_2}= \bar{l_2}\bar{l_1}$ (recall that for any
unit quaternion $q$, $\bar q = q^{-1}$).
A number of other properties can easily be obtained.
If $l=r$, then the real ($x_1$) axis is preserved, so this special
case represents an element of $\SO{3}$ acting on
$(x_2,x_3,x_4)$. 

Elements in $\SO{4}$ are 
either single rotations, that fix all points on a two-dimensional
plane, or double rotations that fix only the origin.
These two types of rotation can be easily be distinguished
in the quaternion notation (see Lemma~\ref{lem:fixedpointspace_order_two_2}).

Using the map  (\ref{quaternionmap}), Conway and Smith \cite{CS} use
the classification of subgroups
of $\OO{3}$ to present a complete list of subgroups of $\SO{4}$.

\section{Series of Groups}\label{sec:series}
In a search for examples of groups with an absolutely irreducible action on a
finite dimensional space, where all the isotropy subgroups have even
dimensional fixed point spaces, we came
across three groups of order $48$ having this property. 
Further investigation showed that these formed part of three infinite
series of such groups. 

We are interested in three series of  groups $G_j(m), j=1,2,3$ and $m\ge 3$ an
odd integer. We follow the notation
in \cite{CS}, where each group is related to two
subgroups of $\OO{3}$. The groups can also be defined by a
(non-minimal) set of generators in the quaternion notation described
above (see \cite{CS}, Tables 4.1 and 4.2).
Let us write
\begin{eqnarray*}
  G_1(m)&=&\pm \frac12[D_{2m}\times \overline{D}_8]\\ 
  G_2(m)&=&\pm \frac14[D_{4m}\times\overline{D}_8]\\ 
  G_3(m)&=&\pm [D_{2m}\times D_4].
\end{eqnarray*}
The orders of these groups are (see \cite{CS}): $|G_1(m)|=\frac12\cdot2\cdot 2m\cdot
8=16m$, $|G_2(m)|=2\cdot\frac14\cdot 8\cdot 4m=16m$ and
$|G_3(m)|=2\cdot 2m\cdot 4=16m$. So we get group orders
$48,80,112,\dots$, all of which have the form $16+32\cdot \ell$,
$\ell\in\N$ and $m=2\ell+1$. Observe that in the notation of \cite{CS}
the group $D_{2n}$ has $2n$ elements. 
Table \ref{tab:SmallGroupLibrary} translates this notation for small values of $m$ into the
SmallGroupLibrary notation of GAP \cite{GAP}. For some of the
computations this program and its library are extremely useful (some
computational results are given in section~\ref{sec:gap}).
\begin{table}[h]
  \centering
  \begin{tabular}{|l|c|c|c|}
\hline
$m$&$G_1(m)$&$G_2(m)$&$G_3(m)$\\
\hline
3&[48:17]&[48:15]&[48:41]\\
\hline
5&[80:17]&[80:15]&[80:42]\\
\hline
7&[112:16]&[112:14]&[112:34]\\
\hline
9&[144:18]&[144:16]&[144:44]\\
\hline
11&[176:16]&[176:14]&[176:34]\\
\hline
13&[208:17]&[208:15]&[208:42]\\
\hline
  \end{tabular}
  \caption{The SmallGroupLibrary names for our groups for small values of $m$.}
  \label{tab:SmallGroupLibrary}
\end{table}

The main results concerning these groups are collected in the following
theorems. 
\begin{thm}\label{series.thm:properties}
  \begin{enumerate}
  \item Given any two groups  within the same series $G_j(m)$ and
    $G_j(m')$ then, if $m$ divides 
    $m'$, we have 
\[G_j(m)\subset G_j(m').\] 
\item The closure of the union of all groups within one family is a
  compact, one-dimensional Lie-Group
  $G_j$. The groups $G_1$ and $G_2$ are isomorphic. 
\item The infinitesimal generator of $G_j$ is given by 
\[L_j=[i,0].\]
\item Each of the groups $G_j(m)$, $j=1,2,3$  contains a unique index
  $2$ subgroup 
  $F_j(m)$ which  commutes with 
\[J=\left(
  \begin{array}{cccc}
   0&1&0&0\\
   -1&0&0&0\\
   0&0&0&1\\
  0&0&-1&0 
  \end{array}\right) = [i,1]
.\]
\item The closure of the unions of all $F_j(m)$  are again
  subgroups of $\SO{4}$, denoted by $F_j$, which are compact one
  dimensional Lie Groups which are contained in $G_j$ and which
  commute with $J$.
\item The elements in $G_j(m)\setminus F_j(m)$, and in $G_j\setminus
  F_j$ anti-commute with $J$. 
  \end{enumerate}
\end{thm}
The following theorem describes the actions of these groups on $\R^4$.
\begin{thm}\label{series.thm:action}
  \begin{enumerate}
  \item The natural actions $\rho$ of $G_j(m)$ on $\R^4$ for $j=1,2,3$
    and $m\ge 3$ are absolutely irreducible. 
  \item If $m\ge 3$ is odd then corresponding to the natural
    representation $\rho$ of $G_j(m)$ 
    on $\R^4$ there exists 
    at least one nontrivial isotropy subgroup. All isotropy subgroups
    are of order $2$ and the corresponding fixed point
    space is $2$-dimensional.
  \item In each of the groups $G_j(m)$ we have precisely $j$ isotropy types.
  \item The normalizer of the isotropy subgroups acts  on
    the fixed point subspaces in the following way:
    \begin{enumerate}
    \item $j=1$: The normalizer is isomorphic to $\DD{2}$ and it acts
      on $\Fix{\ZZ{2}}$ as $\ZZ{2}$, namely as a rotation by $\pi$.
    \item $j=2$: here we have two isotropy subgroups: in one case the
      normalizer is isomorphic to $\DD{2}$ and it acts as in the
      previous case. The other normalizer is isomorphic to $\ZZ{2m}\times\ZZ{2}$
      and it acts on $\Fix{\ZZ{2}}$ as a rotation by $\frac{\pi}{m}$.
    \item $j=3$: in this case we have three isotropy types; each one is
      represented by a group isomorphic to $\ZZ{2}$.  In
      each case the normalizer is isomorphic to $\ZZ{4}\times \ZZ{2}$ and the
      normalizer acts as $\ZZ{4}$.
    \end{enumerate}
  \end{enumerate}
\end{thm}

\section{Flows}\label{sec:flow}
In this section we want to look at the set of $G_j(m)$-equivariant
vector fields on $\R^4$ and we show that generically we find that loss
of stability leads to bifurcating equilibria.   In order to determine
the fine structure of the equivariant maps we need to know the number
of equivariant polynomial maps in a given dimension. Computing the
Poincare series gives this information. The dimension of the space of
equivariant polynomials can be computed directly using a formula given
in Sattinger \cite{Sat1}.
The computations needed 
for such a detailed study are given in Section \ref{characters}.  
Tables \ref{tab:equiv_info1}, \ref{tab:equiv_info2} and
\ref{tab:equiv_info3} show the number of invariant polynomials and  
equivariant polynomial maps respectively in the various degrees for
the groups $G_j(m), F_j(m)$ for $j=1,2,3$ and $m\ge 3$, $m$ odd.
\begin{table}[ht]
  \centering
  \begin{tabular}[h]{|cccccc|ccccccc|}
\hline
 m& $G_1(m)$&e3&i4&i6&i8&$F_1(m)$&e1&i2&e3&i4&i6&i8\\
\hline
 3&[48,17]&3&2&4&9&[24,11]&2&1&6&3&6&15\\
 5&[80,17]&3&2&3&5&[40,11]&2&1&6&3&4&7\\ 
 7&[112,16]&3&2&3&5&[56,10]&2&1&6&3&4&7\\
 9&[144,18]&3&2&3&5&[72,11]&2&1&6&3&4&7\\ 
 11&[176,16]&3&2&3&5&[88,10]&2&1&6&3&4&7\\ 
 13&[208,17]&3&2&3&5&[104,11]&2&1&6&3&4&7\\ 
 15&[240,78]&3&2&3&5&[120,33]&2&1&6&3&4&7\\ 
 17&[272,17]&3&2&3&5&[136,11]&2&1&6&3&4&7\\ 
 19&[304,16]&3&2&3&5&[152,10]&2&1&6&3&4&7\\ 
21& [336,103]&3&2&3&5&[168,41]&2&1&6&3&4&7\\ 
\hline
  \end{tabular}
  \caption{The information on the invariants/equivariants for the
    groups in $G_1(m)$. Here $e$ stands for equivariants, $i$ for
    invariants and the number behind these letters for the degree of
    the polynomial map. The number in the table gives the dimension of
    the space of equivariants/invariants in the given
    degrees. Observe, here and in the following tables the groups in
    the left column
  act absolutely irreducibly and hence we always have $e1=i2=1$.} 
  \label{tab:equiv_info1}
\end{table}
\begin{table}[ht]
  \centering
  \begin{tabular}[h]{|cccccc|ccccccc|}
\hline
 m& $G_2(m)$&e3&i4&i6&i8&$F_2(m)$&e1&i2&e3&i4&i6&i8\\
\hline
 3&[48,15]&3&2&5&9&[24,10]&2&1&6&3&8&15\\
 5&[80,15]&3&2&3&5&[40,10]&2&1&6&3&4&7\\ 
 7&[112,14]&3&2&3&5&[56,9]&2&1&6&3&4&7\\
 9&[144,16]&3&2&3&5&[72,10]&2&1&6&3&4&7\\ 
 11&[176,14]&3&2&3&5&[88,9]&2&1&6&3&4&7\\ 
 13&[208,15]&3&2&3&5&[104,10]&2&1&6&3&4&7\\ 
 15&[240,76]&3&2&3&5&[120,32]&2&1&6&3&4&7\\ 
 17&[272,15]&3&2&3&5&[136,10]&2&1&6&3&4&7\\ 
 19&[304,14]&3&2&3&5&[152,9]&2&1&6&3&4&7\\ 
21&[336,101]&3&2&3&5&[168,40]&2&1&6&3&4&7\\ 
\hline
  \end{tabular}
  \caption{The information on the invariants/equivariants for the groups in $G_2(m)$. Here $e$ stands for equivariants, $i$ for invariants and the number behind these letters for the degree of the polynomial map. The number in the table gives the dimension of the space of equivariants/invariants in the given degrees.}
  \label{tab:equiv_info2}
\end{table}
\begin{table}[ht]
  \centering
  \begin{tabular}[h]{|cccccc|ccccccc|}
\hline
 m& $G_3(m)$&e3&i4&i6&i8&$F_3(m)$&e1&i2&e3&i4&i6&i8\\
\hline
 3&[48,41]&3&3&4&10&[24,10]&2&1&6&3&8&15\\
 5&[80,42]&3&3&3&6&[40,10]&2&1&6&3&4&7\\ 
 7&[112,34]&3&3&3&6&[56,9]&2&1&6&3&4&7\\
 9&[144,44]&3&3&3&6&[72,10]&2&1&6&3&4&7\\ 
 11&[176,34]&3&3&3&6&[88,9]&2&1&6&3&4&7\\ 
 13&[208,42]&3&3&3&6&[104,10]&2&1&6&3&4&7\\ 
 15&[240,182]&3&3&3&6&[120,32]&2&1&6&3&4&7\\ 
 17&[272,43]&3&3&3&6&[136,10]&2&1&6&3&4&7\\ 
 19&[304,34]&3&3&3&6&[152,9]&2&1&6&3&4&7\\ 
21&[336,201]&3&3&3&6&[168,40]&2&1&6&3&4&7\\ 
\hline
  \end{tabular}
  \caption{The information on the invariants/equivariants for the groups in $G_3(m)$. Here $e$ stands for equivariants, $i$ for invariants and the number after these letters for the degree of the polynomial map. The number in the table gives the dimension of the space of equivariants/invariants in the given degrees.}
  \label{tab:equiv_info3}
\end{table}
Looking at the lowest order nontrivial $G_j(m)$-equivariant polynomial
maps, we observe that we expect in each case three independent maps of
order three. In the cases $j=1,2$ we expect two of these maps to be
variational. This is because $i4=2$ and the gradient of an invariant
of degree 4 is an equivariant of degree 3 that is variational.  In the
case $j=3$ all of the equivariants are variational. For fixed 
$j$ the numbers $e3$, $i4$, $e5$ (which are not displayed in the
tables) are monotonically decreasing in $m$. So if we find $3$
independent equivariant maps of order 3, which are equivariant for all
$m$, we see that these three maps are the ones to be looked
at. Moreover they are equivariant with respect to $G_j$, which is a
compact Lie group. In order to discuss specifics for each group, one
has to look at higher order equivariants. However we expect, that
generically the bifurcation scenario will be decided at the cubic level. In
the case $j=3$, the equivariant maps up to order $3$ are variational,
so restricting to third order we will have bifurcation to equilibria. 
We collect the results in the following two theorems, which will proved in
Section \ref{proofs}.
\begin{thm}\label{thm:equivariant_structure}
  For each $j=1,2,3$ and each $m\ge 3$, $m$ odd, there are precisely
  three linearly independent cubic equivariant maps. For $j=1,2$ two
  of these maps are gradients of invariant polynomials, the third one
  is a Hamiltonian vectorfield. In the case $j=3$ all three
  vector fields are gradients of invariant polynomials.  
\end{thm}
\begin{thm}\label{thm:bifurcation}
  The third order polynomial equations lead to bifurcation to one or
  more circles 
  of equilibria. At least one of these circles intersects one fixed
  point spaces in 
  discrete points. Each of these points is a regular point, so if we
  restrict the map to the complement of a ball around zero in the
  nontrivial fixed point spaces, we find at least one point  which
  persists under perturbation with higher order terms. 
\end{thm}
In a short form we have shown:
\begin{thm}\label{thm:persistance}
  Generically, i.e. for an open and dense set in the set of
  $C^\infty$-equivariant vector fields we have bifurcation of nontrivial
  equilibria. 
\end{thm}
In this sense the Ize conjecture holds  for the groups under consideration.
\section{Hamiltonian structure}\label{sec:Hamilton}
In this section we want to describe a Hamiltonian structure which we
have in all the groups discussed here, but has not been well
studied. Whether the generic behaviour for this type of group is
different from the usual one is not clear, it could well be that there
are new phenomena. Let us briefly describe the situation of
equivariant Hamiltonians.
On $\R^{2n}$ we look at linear operators $J$ with $J^2=-\1$. We call a
vectorfield $v:\R^{2n}\to\R^{2n}$ Hamiltonian, if there exists a
function $H:\R^{2n}\to\R$ with
\[v=J\nabla H.\]
The vectorfield $v$ is equivariant with respect to a group $G$ if one
of two conditions hold:
\begin{enumerate}
\item $H$ is an invariant for $G$, and $J$ commutes with $g\in
  G$. Then $\nabla H$ is equivariant, and $v(gx)=J\nabla
  H(gx)=Jg\nabla H(x)=gJ\nabla H(x)=gv(x)$. Observe that $J$ is not a
  multiple of the identity and commutes with $g\in G$, therefore the
  action of $G$ is not absolutely irreducible.
\item In this case we require an index $2$ subgroup $F$ of $G$ and $H$
  is an invariant of $F$, $J$ commutes with $F$ and for $g\in
  G\setminus F$ we have $H(gx)=-H(x)$ and $g\inv Jg =-J$. Then,
  obviously $v$ is equivariant. In such a case $G$ can act absolutely
  irreducibly.   
\end{enumerate}
In all our examples we have a pair $(G,F)$ of index $2$ subgroups
and we have functions $H$ and operators $J$ which are invariant under
$F$ and anti-commute with elements in $G\setminus F$. The details can
be found in the next section. 
 \section{Proofs}\label{proofs}
\cite{CS} give a set of generators for these groups. 
For elements $a,b,c,\dots$ of a group we write
\[\langle a,b,c,\dots\rangle\]
for the smallest subgroup containing these elements.
Following \cite{CS} we define elements using the short notation
\[e_s=e^{\frac{\pi\,i}s}.\]
With this notation the groups are given by (see \cite{CS}, Tables 4.1
and 4.2)
\begin{eqnarray*}
G_1(m)&=&\langle [e_m,1],[1,i],[1,j],[j,e_{4}]\rangle\\
G_2(m)&=&\langle [e_m,1],[1,i],[e_{2m},j],[j,e_4]\rangle\\
G_3(m)&=&\langle [e_m,1],[1,i],[j,1],[1,j]\rangle. 
\end{eqnarray*}
Observe that these sets of generators do not form  minimal sets of
generators. A 
first observation is the following: if we multiply the first two
elements in $Q\times Q$, then the product generates the same group as
the two elements. For this it suffices to prove, that the element
$[e_m,1]$ is in the group generated by $[e_m,i]$. Since 
\[[e_m,i]^4=[e_m^4,1]\]
and $m$ and $4$ are relatively prime, $e_m$ and $e_m^4$ generate the
same group in $S^1$. 
\begin{lem}\label{lem:cyclic_group}
  The element $[e_m,i]$ generates a group $H(m)$ of order $4m$. 
\end{lem}
\pf The order of the group generated by this element is obviously a
multiple of $m$. So we have
\[[e_m,i]^m=[-1,\pm i]\not=[1,1], \ \ [e_m,i]^{2m}=[1,(-1)^m]\not=[1,1],\mbox{
  and }[e_m,i]^{4m}=[1,1].\]
\endproof~\\[3mm]

\subsection{Proof of Theorem \ref{series.thm:properties}.}

In this section we give the proof of the six parts of
Theorem~\ref{series.thm:properties} on the structure of the groups.

\begin{enumerate}
\item 
If $m$ divides $m'$, so $m'=pm$ for some integer $p$, then 
\[
[e_{m'},1]^p = [e^{p\pi i/m'},1] = [e^{\pi i/m},1] = [e_{m},1]
\]
so $\langle [e_{m},1] \rangle \subset \langle [e_{m'},1]\rangle$
and hence $G_1(m) \subset G_1(m')$ and  $G_3(m) \subset G_3(m')$.
Similarly, considering the element  $[e_{2m},j] \in G_2(m)$, 
$[e_{2m'},j]^p = [e_{2m},j^p]$, but we know that $p$ is odd (since $m'$
is odd), so this is  $\pm [e_{2m}, j] \in G_2(m)$, 
so $\langle [e_{2m},j] \rangle \subset \langle [e_{2m'},j]\rangle$
and hence $G_2(m) \subset G_2(m')$. 
\item The closure of the union of all the groups $G_1(m)$ is 
\[
G_1 = \langle [e^{i\theta},1],[1,i],[1,j],[j,e_{4}]\rangle
\]
where $\theta \in [0,2\pi]$. This is a compact one-dimensional Lie
group. 
For the group $G_2$, we also have elements of the form 
$[e^{i\phi},j]$, but these are already included in $G_1$,
so $G_2\subset G_1$. But also  $G_1\subset G_2$, since $[1,j]\in G_2$,
so $G_1 = G_2$. 
By the same argument the closure of $G_3(m)$ is a compact Lie group,
but this is not the same group as $G_1$. 
\item 
Note that $H(m)$ is in the intersection of all $G_j(m)$. The closure
of the union
over all groups $H(m)$ is a one-parameter group and hence isomorphic
to $S^1$ with the generator of its Lie algebra given by $[i,0]$. The
group generated by $H(m)$ and the remaining generators of $G_j(m)$
produce a extension of finite index, therefore the connected component
of this group is isomorphic to $S^1$.
\item 
We define
\begin{eqnarray*}
F_1(m)&=&\langle [e_m,i],[1,j]\rangle\\
F_2(m)&=&\langle [e_m,i],[e_{2m},j]\rangle\\
F_3(m)&=&F_1(m).
\end{eqnarray*}
Let us write
\[F(m)=\langle [e_m,i]^2\rangle = \langle [e_{m/2},-1]\rangle.\]
Clearly $F(m)\subset F_j(m)$. The generator of $F(m)$ commutes with
all elements in $F_j(m)$ and hence $F(m)$ is contained in the center
of $F_j(m)$ for $j=1,2,3$. $F(m)$ contains $2m$ elements, including
minus the identity. 
Since $[e_m,i]$ does not commute with the second
generator of $F_j(m)$, $F(m)$ is the center of $F_j(m)$. 
Now the square of the second generator  
of $F_j(m)$ is in both cases in $F(m)$, since for $F_1(m)$ the square of
the second generator is $-\1 \in F(m)$ and for $F_2(m)$ the square of
the second generator is $[e_m,-1] = \pm [e_{m/2},-1]^{(m+1)/2} \in
F(m)$. 
Therefore  $F_j(m)$ is an index $4$ extension of $F(m)$ and hence
it has $8m$  elements. Therefore $F_j(m)$ is index $2$ subgroup of
$G_j(m)$.
\item It is obvious that each element of $F_j(m)$ commutes
with $[i,1]$. Therefore $F_j(m)$ commutes with $J$ which is the map
induced by $[i,1]$.
\item The elements of $G_j(m)$ which are not in $F_j(m)$
anti-commute with $[i,1]$.  Note that $F(m)$ being the center of
$F_j(m)$ is normal in $G_j(m)$ and $|G_j(m)/F(m)|=8$.
\end{enumerate}
\endproof
Before we enter the proof of Theorem \ref{series.thm:action}, we state
some useful little lemmas which should be well known, but we could not find
a reference. 
\begin{lem}\label{abs_irred_orth_crit}
  If a Lie Group $G$ acts on a real space $V$  and the condition 
\[(\forall A\in\OO{V}\;\;:\;\;Ag=gA)\;\; \imp\;\; A=\pm\1\]
is true, then the action is absolutely irreducible.
\end{lem}
\pf First observe that the action is irreducible: assume $U\subset V$ is a
$G$-invariant subspace with an orthogonal complement
$W$. Then the orthogonal projection onto $U$ along $Q$ and vice versa commute
with $G$. Especially the operator $I$ which acts as $\1$ on $U$ and as $-\1$ on
$W$ commutes with $G$. But this operator is in $\OO{V}$ and therefore  $I=\pm
\1$ and one of these spaces is $\{0\}$ and the other equal to $V$. \\
Now the set of commuting matrices forms a division algebra and if the action
is not absolutely irreducible it contains a subspace isomorphic to $\C$. let
$J$ be the operator corresponding to $i$, then $J$ is skew and $J^TJ=\1$. Let
$\alpha, \beta \in \R$, $\alpha\not=0,\ \beta\not=0$ with
$\alpha^2+\beta^2=1$. Then $\alpha\1+\beta J$ 
commutes with $G$ and it is orthogonal, since
\[(\alpha\1+\beta J)^T(\alpha \1+\beta J)=
\alpha^2\1+\alpha\beta(J^T+J)+\beta^2J^TJ=(\alpha^2+\beta^2)\1=\1.\]
Therefore $\alpha\1+\beta J$ is a multiply of the identity and therefore we
have a contradiction.  
\endproof
\begin{lem}\label{lem:order_two}
  Elements $g\in \SO{4}$ with $g^2=\1$ and which are not equal to $\pm
  \1$ have a two-dimensional fixed point space. 
\end{lem}
\pf All eigenvalues $\lambda$ of $g$ satisfy $\lambda^2=1$ and hence
they are equal to $\pm 1$. $\det g=1$ implies that the number of
eigenvalues equal to $-1$ is even, and so this number is $0,2,4$. The
cases $0,4$ are excluded by our assumption $g\not=\pm\1$.
\endproof
\begin{lem}\label{lem:fixedpointspace_order_two}
  Let $a,b$ be two unitary quaternions, such that $[a,b]$ has order two
  and $[a,b]\not=\pm\1$. Then $[a,b]$ fixes the two elements $1+a\inv
  b$, $a+b$.  
  These two elements span a two dimensional subspace.
\end{lem}
\pf $[a,b]^2=1$ means $[a^2,b^2]=[1,1]$ or $[-1,-1]$. In either case,
$a^2=b^2$.  We observe that 
\[[a,b](1+a\inv b)=a\inv b +a^{-2}b^2=a\inv b+1.\] 
In the same fashion we get
\[[a,b](a+b)=b+a\inv b^2=b+a.\] 
This leads to a  two-dimensional space except in the case 
\[(a+b)=r(1+a\inv b)\]
for some real number $r\in\R$. Since the nonzero quaternions form a
multiplicative group, we have $r=a$. Therefore $a\in\R$. Since
$[a,b]^2=\1$, we have $a^2=\pm 1$ and so $a=\pm 1$ and
$a^2=1$. Since $b^2=a^2=1$,  $b=\pm 1$. Then we have
\[ [a,b]=[\pm 1,\pm 1]=\pm\1\]
which contradicts our assumption.
\endproof
\begin{lem}\label{lem:fixedpointspace_order_two_2}
  An element $[l,r]$, $l,r\in Q$ fixes an element $p$ if and only if $l$ and
  $r$ are conjugate.
\end{lem}
\pf If $[l,r]p=p$ then 
\[l\inv p r=p\]
or
\[r=p\inv l p.\]
So $l$ and $r$ are conjugate.
Note that $p $ is mapped to $p$
and $l p$ is mapped to $l p$, so we have a two-dimensional fixed-point
space (except in the case $l=1$, but in that case $[l,r]$ is the
identity).
\endproof
This lemma is useful to determine whether elements have fixed point
subspaces. It also distinguishes between the single rotations and
double rotations in $\SO{4}$.  To make use of this we need the following observation.
\begin{lem}\label{lem:conjugacy}
  Two elements in $Q$ are conjugate if and only if they have the same
  real part.
\end{lem}
\pf This follows from the work of Janowsk\`a and Opfer
\cite{JaOp}. They prove that two quaternions are conjugate, if they
have the same length and same real parts.
\endproof

It is also  useful to have an explicit form of the space fixed by
a single rotation. 
\begin{lem}\label{lem:rotfix}
Let $[l,r]$ be a single rotation, with
$Re(l)=Re(r)$.  Then $[l,r]$ fixes the space spanned by the two
vectors $l - \bar r$ and $1 - \bar l \bar r$
\end{lem}
\pf
$[l,r]$ maps $l - \bar r$ to $\bar l (l - \bar r) r = r-\bar l$.
But if $Re(l)=Re(r)$, then $l+\bar l = r + \bar r $, 
so $l - \bar r =  r-\bar l$.
Similarly, $1 - \bar l \bar r$ is mapped to 
$\bar l (1 - \bar l \bar r) r = \bar l r - \bar{l}^2$
but since $l+\bar l = r + \bar r $, we have 
$1 + \bar{l}^2 = \bar l r + \bar l \bar r$ so this mapping is also the
identity. 
\endproof

For example, consider the element $[i,j]$, which is of order 2. 
The real part of both quaternions is zero so this is a single
rotation. The fixed vectors are $i+j$ and $1-k$.

\subsection{Proof of Theorem~\ref{series.thm:action}}

1. To prove that the action of $G_j(m)$ on $\R^4$ is absolutely
  irreducible, we use  Lemma~\ref{abs_irred_orth_crit}.
Each of our groups  $G_j(m)$ contains the element $[1,i]$.
 If $[l,r]\in Q\times Q$ commutes with  $[1,i]$ then 
$r=q_1+q_2 i$. Now we also have an element of the form $[*,j]$ in the
group, and  $r\cdot j=j\cdot r$ implies $q_2=0$ showing that the right
element $r$ is real.\\
Now we prove a similar statement for $l$. Each group contains
  $[e_m,1] = [\cos(\pi/m)+i\sin(\pi/m),1]$. Since $m\ge 3$,  the property
$e_m\cdot l=l\cdot e_m$
implies that  $l=p_1+ip_2$. Now we also have in each group an element
of the form $[j,*]$, implying $p_2=0$. Therefore the only commuting elements
are of the form $[p_1,q_1]$ with $p_1,q_1\in\R$. Since $l$ and $r$ are
unit quaternions, we conclude $p_1,q_1=\pm 1$. Therefore all the
commuting elements in $\OO{4}$ are of the form $\pm\1$ and we deduce
absolute irreducibility from Lemma~\ref{abs_irred_orth_crit}.

To prove our main results concerning the isotropy subgroups of
$G_j(m)$,  we begin with some general observations which are relevant
for all three groups; the second part of the proof will address each
group separately. \\

Note first that none of the nontrivial elements of the group
$H(m)=\langle[e_m,i]\rangle$ of order $4m$ fixes 
any point. Using Lemma~\ref{lem:fixedpointspace_order_two_2}, $[e_m,i]^r$ can only fix
a subspace if $Re(e_m^r)$ = $Re(i^r)$, which implies that 
\[\cos(\pi r/m)=Re(i^r).\]
If $r$ is odd this equation cannot be satisfied since the right-hand
side is zero but the left-hand side is not zero, because $m$ is odd.
If $r$ is even the condition can only be satisfied if $r$ is a
multiple of $m$. If $r=2m$ then the left-hand side is 1 and the
right-hand side is $-1$. For  $r=4m$ the equation is satisfied but
that element is the identity. For any isotropy subgroup we conclude
that it intersects $H(m)$ only at the trivial element.

The remaining argument for parts 2,3 and 4 of the theorem is different
for the three groups in question and we discuss each case separately. \\
(a) The case $G_1(m)$:\\
2. We consider the subgroup $F_1(m)$ generated by $H(m)$ and the element
$[1,j]$. It was shown in the proof of
Theorem~\ref{series.thm:properties} that $F_1(m)$ is of order
$8m$. Hence $H(m)$ is a normal subgroup of $F_1(m)$ and therefore any isotropy
subgroup in $F_1(m)$ has order $2$ and the nontrivial element is in
the coset  $H(m)\cdot[1,j]$. The elements in this coset have the form
\[[e_m^r,i^rj].\]
None of these elements has order $2$, because $i^r j = \pm j$ or $\pm
k$ which when squared gives $-1$, but $e_m^{2r} \ne -1$, and
therefore, none of these 
elements fixes any nontrivial $x$. So, isotropy subgroups are
subgroups of $G_1(m)$ and intersect $F_1(m)$ only at the identity. Hence
all isotropy subgroups have order $2$, and by
Lemma~\ref{lem:order_two} have a two-dimensional fixed point space. 
The nontrivial element lies in
  the nontrivial coset of $F_1(m)$. So it has the form
\[[e_m^r j,i^r e_4] \mbox{ or }[e_m^r j,i^r j e_4].\]
Since the square of the second component of the first element is never
$\pm 1$ the first element is never of order $2$. So, we concentrate on
the second element. We begin with two simple remarks:
\[j1j=-1, jij=i,\]
and therefore for $z\in \C$ we find
\[jzj=-\bar z.\]
From this we conclude for $|z|=1$ that $zjzj=-1$. 
In particular, for $p\in\N$ we get
\[e_pje_pj=-1.\]
The square of the second element has the form
\[[e_m^rje_m^rj,i^rje_4i^rje_4]=[-1,i^rji^re_4je_4]=[-1,je_4je_4]=[-1,-1]=
[1,1].\] 
Therefore this coset of $F_1(m)$ consists of $4m$ elements of order $2$.
By Lemma~\ref{lem:order_two}, each of these elements generates an
isotropy subgroup isomorphic to $\ZZ{2}$ with a two-dimensional fixed
point space. 
The fixed point spaces for these elements are given by Lemma
\ref{lem:fixedpointspace_order_two}:\\
with $a=e_m^r j$, $b=i^r j e_4$ we get the fixed point space is spanned by
\[e_m^r j+i^r j e_4,\ \ \ 1+(e_m^r j)\inv i^r j e_4.\]
3. From 
\begin{eqnarray*}
[e_m,i][j,je_4]&=&[e_mj,ije_4]=[j\overline{e_m},-jie_4]\\
&=&[j\overline{e_m},je_4\bar i]\\
&=&[j,je_4][e_m,i]^{-1}  
\end{eqnarray*}
we get
\[[e_m,i]^2[j,je_4]=[e_m,i][j,je_4][e_m,i]^{-1}.\]
This shows that at least $2m$ of these $4m$ elements are conjugate
under $F(m)$. 
Consider
\[[1,-j][j,je_4][1,j]=[j,-jje_4j]=[j,e_4j]\]
which is another conjugate element of order $2$. Observe that for
$q=1$ or $q=3$ 
\[[e_m^{qm}j,i^{qm}je_4][j,e_4j]=[-j^2,i^{qm}je_4e_4j]=[1,i^{qm}jij]=[1,i^{qm+1}]=[1, 
1]\]
where $q=3$ if $m=1\mod 4$ and $q=1$ if $m=3\mod 4$. 
This proves that all elements of order $2$ in this coset are conjugate
and hence there is precisely one isotropy type with two dimensional
fixed point space. 

4. Let $S_1(m)$ be a representative of this isotropy type.
Now we have $4m$ objects (either the subgroups or their invariant
planes) that are permuted by the
group $G_1(m)$, so by the orbit-stabilizer theorem, the 
stabilizer of any of these objects must be of order $4$. 
The stabilizer is also the normalizer of $S_1(m)$, that is, the
largest subgroup of $G_1(m)$ in which  $S_1(m)$ is normal.
Clearly the stabilizer includes $-\1$, so the stabilizer is isomorphic
to $\DD{2}$ and acts on $\Fix{S_1(m)}$ as minus the identity, i.e.\ as
a rotation through $\pi$.

(b) The case $G_2(m)$:\\
 Here we follow the lines of the previous argument. None of the
 elements in $H(m)$ fixes anything. $H(m)$ is an index $2$ subgroup
 of $F_2(m)$ and therefore every isotropy subgroup of $F_2(m)$ is of
 order $2$ and intersects $H(m)$ in the trivial element. Again we
 look for order $2$ elements in $F_2(m)$. The nontrivial coset of
 $H(m)$ in $F_2(m)$ is the coset of $[e_{2m},j]$, the general element
 in this coset therefore is given by
\[[e_m^r,i^r][e_{2m},j],\ \ r=0,\dots,4m-1.\]
Squaring these elements gives us
\[[e_{2m}^{2r}e_{2m},i^rj]^2=[e_{2m}^{4r+2},i^rji^rj]=[e_{2m}^{4r+2},-1].\]
So we are interested in those $r$ with 
\[e_{2m}^{4r+2}=-1\]
or
\[4r+2=2m\mod 4m\iff 2r+1=m\mod 2m. \]
We get four solutions in the set $0\le r\le 4m-1$: writing $m=2\tau+1$
then (obviously) the solutions $r$ have the form $r=\tau\mod m$.
\[r_j=\tau+q m,\ \mbox{ for }q=0,1,2,3.\]
From Lemma \ref{lem:order_two} it follows that the dimension of the
corresponding fixed point space is two. 
Depending on the parities of $q$ and $\tau$
the exponent $\tau+qm$ can be even or odd: for each parity of $\tau$
there are two parities of $q$ leading to an odd exponent and also two
parities leading to an even exponent. In the odd case  the element
takes the form
\[[i,\pm k],\]
in the even case the form
\[[i,\pm j].\]
In any case the elements 
\[[i,k]\mbox{ and }[i,-k]\]
and the elements
\[[i,j]\mbox{ and }[i,-j]\]
are conjugate in $F_2(m)$ (the conjugating element is $[1,i]$). 
The elements $[i,j]$ and $[i,k]$ are not
conjugate (as one can easily check) in $F_2(m)$,
so we get two isotropy types in $F_2(m)$. 
However these two isotropy subgroups in $F_2(m)$ are conjugate within
$G_2(m)$, since $[j,ie_4][i,j][-j,-i\overline{e_4}]=[i,k]$. Hence we
have so far only found one isotropy type in $G_2(m)$. 
\\
Now we have to look at the full group. $F_2(m)$ is an index $2$
subgroup of $G_2(m)$ and therefore any isotropy subgroup of $G_2(m)$
intersects $F_2(m)$ in the trivial group or one of the isotropy
subgroups of order $2$ in $F_2(m)$. Therefore we have to look for
elements of order $2$ and $4$ in the coset of $F_2(m)$. This coset
consists of $8m$ elements $f[j,e_4]$ where $f\in F_2(m)$. So we get
this coset as a union of two sets
\[\Meng{[e_m^r,i^r][j,e_4]}{0\le r<4m}\]
and 
\[\Meng{[e_m^re_{2m},i^rj][j,e_4]}{0\le r<4m}.\]
So the first class of elements has the form
\[[e_m^rj,i^re_4]\]
and the squares are of the form
\[[e_m^rje_m^rj,i^{2r}i]=[-1,\pm i]\not=[1,1].\]
Squaring again gives  
\[[1,-1]\not=[1,1]\]
and so there are no elements of order $2$ or $4$ within this class,
therefore none of these elements belongs to any isotropy subgroup. \\
The second class of elements has the form
\[[e_{2m}^{2r+1}j,i^rje_4].\]
The squares have the form
\[[ e_{2m}^{2r+1}je_{2m}^{2r+1}j,i^rje_4i^rje_4]=[-1,(-1)^ri^rje_4je_4i^r].\]
This gives 
\[[-1,(-1)^r(-1)(i^r)^2]=[-1,-1]=[1,1].\]
None of these elements of order $2$ can be conjugate to any of the 
four in $F_2(m)$ we had
found before, since $F_2(m)$ is a normal subgroup of $G_2(m)$. \\
To show that they are conjugate to each other, consider 
\[ [e_m,i][e_{2m}j,je_4][\overline{e_m},-i]= [e_m,i]^2[e_{2m}j,je_4].\]
As before this type of conjugation gives us two classes
\[ [e_m,i]^{2r} [e_{2m},j][j,e_4]\]
and \[ [e_m,i]^{2r+1} [e_{2m},j][j,e_4].\]
Conjugation with $[e_2m,j]$ shows that the elements in these two classes
are all conjugate, since
\[
[e_{2m},j] [e_{2m},j][j,e_4] [\overline{e_{2m}},-j]=[e_{2m}^3 j, e_4 j]
 = [e_m,i][e_{2m},j][j,e_4].
\]
So in total we find in $G_2(m)$ two classes of order $2$ subgroups. 
By Lemma \ref{lem:order_two} the fixed point space is two dimensional
and is given according to Lemma \ref{lem:fixedpointspace_order_two}. 

$G_2(m)$ acts as a permutation group on these elements. The first
class has length $4$, and so $4m$ elements fix each of the groups
isomorphic to $\ZZ{2}$, i.e.\ the normalizer of each subgroup is of
order $4m$. Consider the representative $\langle [i,j]\rangle$ of this
class. Since $[e_{2m},j]$ and $-\1$ commute with $[i,j]$, and
$[e_{2m},j]$ is of order $2m$, these two elements generate a group
$\ZZ{2m}\times \ZZ{2}$ of order $4m$ which must be the normalizer. The
action of the normalizer on the two-dimensional space is the
normalizer quotient which acts as $\ZZ{2m}$, a rotation through an
angle $\pi/m$.  \\
In the other class of isotropy subgroups we have $4m$ elements, so
each one is fixed by $4$ elements and the normalizer acts as
$\DD{2}/\ZZ{2}=\ZZ{2}$, as in the case of $G_1(m)$.  

(c) The case $G_3(m)$:\\
The argument here is very similar to that for $G_1(m)$. 
It has already been shown that there are no isotropy subgroups in 
$F_3(m)=F_1(m)=\langle [e_m,i],[1,j]\rangle$, which is an index 2
subgroup of $G_3(m)$.
Therefore isotropy subgroups of $G_3(m)$ can
only be of order 2 and must be generated by an element in the
nontrivial coset of $F_3(m)$ in $G_3(m)$, $F_3(m) [j,1]$. 
We find $6m$ elements of order $2$:
\[[e_{m}^rj,\eps],\mbox{ where }r\in\{0,\dots,2m-1\}\mbox{ and
}\eps\in\{i,j,k\}.\] Each of these elements is clearly of order $2$
and hence generates a subgroup isomorphic to $\ZZ{2}$.
The second entry determines the conjugacy class,
since there is no possible conjugating element in $G_3(m)$ that
could alter the second entry.
In fact any two of those elements with the same second element are conjugate. 
So the subgroups of order two come in three
conjugacy classes, each of the subgroups has a two-dimensional fixed point
space and the length of each conjugacy class under $G_3(m)$ is $2m$. Each element is fixed by
$8$ elements, the normalizer of the group $\Sigma=\langle
[e_m^rj,\eps]\rangle$ has the form
\[N_{G_3(m)}(\Sigma)=\{[1,1],[1,-1],[e_m^r,\eps],[1,\eps],[e_m^rj,1],
[-1,\eps],[j,-1],[-j,\eps]\}.\]
Therefore it consists of a  group of order $8$ with two generators,
isomorphic to $\ZZ{4}\times \ZZ{2}$. The
quotient is a cyclic group of order $4$. This concludes the proof. 
\endproof 
\nid{\bf Proof of Theorem \ref{thm:equivariant_structure}.}  We begin
to investigate the structure of the invariant polynomials. We will not
give a complete description of all invariants, but just enough to
prove the bifurcation results.
The results
and the arguments are slightly different for the three cases, so we
discuss them partially separately. 
Since the groups $G_j(m)$ operate absolutely irreducibly there is no
linear invariant, this means we look only for invariants which are at
least quadratic.  Let us write
\[I_2(x)=\sum_{\nu=1}^4 x_\nu^2.\]
This is clearly a quadratic invariant for all groups in question. 
Now we restrict our attention to the
groups $F_j(m)$, $j=1,2$. We recall the generating elements:
\[F_1(m)=\langle[e_m,1],[1,i],[1,j]\rangle,\ \
F_2(m)=\langle[e_m,1],[1,i],[e_{2m},j]\rangle.\]  
In both cases we have as one of the generating elements the
element $[e_m,1]$. In order to describe its invariant functions we
introduce complex notation via 
\[z_1=x_1+ix_2,\ \ \ \  z_2=x_3+ix_4 \mbox{ with }
x=x_1+ix_2+jx_3+kx_4=z_1+z_2j.\] 
Then $I_2=|z_1|^2+|z_2|^2$. 
In order to describe further invariants
let us look at the action of $[e_m,1]$ on the complex  variables
$z_1,z_2$. By 
\[[e_m,1]x=\bar e_m x=\bar e_m (z_1+z_2j)=\bar e_m z_1+\bar e_m z_2j.\]
This means the first generator
sends the pair $(z_1,z_2)$ to
$(e^{-\frac{i\pi}m}z_1,e^{-\frac{i\pi}m}z_2)$. The second generator maps 
\[(z_1,z_2)\mapsto(iz_1,-iz_2).\]
For the third one we obtain in the case $F_1(m)$
\[[1,j](z_1,z_2)=(-z_2,z_1)\]
and in the case $F_2(m)$ the third generator maps
\[[e_{2m},j](z_1,z_2)=-(\bar e_{2m}z_2,\bar e_{2m}z_1).\]
We look at monomials in the form $z_1^{k_1}z_2^{k_2}{\bar
  z_1}^{\ell_1}{\bar z_2}^{\ell_2}$. So for the action of the first
element we simply get
\[e^{-\frac{i\pi}m(k_1+k_2-\ell_1-\ell_2)}z_1^{k_1}z_2^{k_2}{\bar
  z_1}^{\ell_1}{\bar z_2}^{\ell_2}.\]
In order that this is invariant under the action of the first element
and which has an order  at most $4$ we have
\[k_1+k_2=\ell_1+\ell_2.\]
This implies that functions of
\[|z_1|^2,|z_2|^2,z_1\bar z_2,\bar z_1z_2\]
are invariant under this particular action.
The second element multiplies the last two expressions with $-1$, so
we should look at the squares of these elements of products of two
sign changing functions, i.e. we look at
functions of
\[|z_1|^2,|z_2|^2,z_1^2\bar z_2^2,\bar z_1^2z_2^2.\] 
Since the third element, in the case of $F_1(m)$, basically
interchanges $z_1, z_2$ the invariants
for $F_1(m)$ have to be symmetric in $z_1,z_2$. Therefore invariants
for $F_1(m)$ have to be functions of
\begin{equation}\label{equ:basic.invariants}
I_2=|z_1|^2+|z_2|^2, I_{4,1}=|z_1|^2|z_2|^2, I_{4,2}=z_1^2\bar
z_2^2+\bar z_1^2 z_2^2, I_6=(|z_1|^2-|z_2|^2)i(z_1^2\bar z_2^2-\bar z_1^2z_2^2).
\end{equation}
In the case of $F_2(m)$ we find it leaves the same functions
invariant, observe that the extra factors multiply to $1$. 

With this information we construct fourth and six order invariants for
$F_j(m)$. Of course $I_2^2$ is a fourth order invariant for $F_j(m)$,
$j=1,2,3$.  Let us consider 
\[I_{4,1}(x_1,\dots,x_4)=|z_1|^2|z_2|^2=(x_1^2+x_2^2)(x_3^2+x_4^2).\]
It is (obviously) invariant  under $F_j(m)$.  Let
\[I_{4,2}=z_1^2\bar z_2^2+\bar z_1^2z_2^2.\]
Let us write
\[I_6=(|z_1|^2-|z_2|^2)i(z_1^2\bar z_2^2-\bar z_1^2z_2^2).\]
This is invariant under $F_j(m)$. We now write down the invariants up
to order $6$ (for $m$ sufficiently large, for small $m$ there might be
additional invariants):
$I_2$ is the unique quadratic invariant, $I_2^2, I_{4,1}, I_{4,2}$ are
the quartic invariants, and
\[I_2^3,I_{4,1}\cdot I_2,I_{4,2}\cdot I_2,I_6\]
are sextic invariants. In a similar way we can construct $7$
invariants of order $8$ and $9$ invariants of order $10$.\\[5mm]
\begin{enumerate}
\item The case $G_1(m)$\\
We provide the explicit form of the generating elements outside
$F_1(m)$ in complex notation.
We get
\[[j,e_4](z_1+z_2j)=-\bar e_4\bar z_1j+e_4\bar z_2.\]
We can write this as $(z_1,z_2)\mapsto (e_4\bar z_2,-\bar e_4\bar
z_1)$. 
 Obviously
$I_2,I_2^2,I_{4,1},I_2^3,I_{4,1}I_2$ are invariant. The function
$I_{4,2}=z_1^2\bar z_2^2+\bar z_1^2z_2^2$ is mapped by the above
substitution to $-I_{4,2}$ (observe $e_4^4=-1$). 
Then $I_{4,2}^2$ is invariant under this
element and hence under $G_1(m)$. 
This substitution applied to $I_6$ changes the signs of both factors
and hence $I_6$ is invariant. 
\item Since the invariants we have constructed are the same for
  $F_1(m)$ and $F_2(m)$ they have to be the same for the groups
  $G_1(m)$ and $G_2(m)$, because the extension is defined with the
  same element. 
\item The case $G_3(m)$. In this case $F_3(m)=F_1(m)$ and therefore we
  just have to look at the remaining generating element $[j,1]$ which
  acts as
\[[j,1](z_1+z_2j)=\bar jz_1+\bar jz_2j=-\bar z_1j +\bar z_2.\]
Therefore we get the substitution $(z_1,z_2)\mapsto (\bar z_2,-\bar
z_1)$. Applying this to the invariants
$I_2,I_2^2,I_{4,1},I_2^3,I_2\cdot I_{4,1}$ we easily see
that they are invariant as well. Note however, that $I_{4,2}$ is
invariant under this element as well, so we find three independent
quartic polynomial invariants in this case.   
\end{enumerate}
If we look at the gradients of $I_2^2,I_{4,1},I_{4,2}$ they give rise
to three independent equivariant maps. In the case of $G_{1}(m)$ and
$G_2(M)$ the 
first two are gradient vector fields, the third gives rise a Hamiltonian
field as we discussed before. In the case $G_3(m)$ we get three
equivariant gradients, so up to cubic level a $G_3(m)$-equivariant
vector fields is a gradient. \\
Observe that our proof so far shows only that there are at least three
equivariant vector fields for the given groups. However, the character
formula for the number of equivariant fields in Section~\ref{characters}
yields three equivariant cubic  fields (compare Tables
\ref{tab:equiv_info1}, \ref{tab:equiv_info2} and
\ref{tab:equiv_info3}). Since for a given degree the number of
equivariant vector fields is as a function of $m$ non increasing, we
conclude that 
we have precisely three cubic equivariants for all $j=1,2,3$ and all
odd $m\ge 3$. 
\endproof~\\[3mm]
{\bf Proof of Theorem \ref{thm:bifurcation}.} We use the standard
theory of complex structures (see e.g. Range \cite{Ran-86}, Chapter
III, Section II) to derive the real vectorfield from the complex form
of the invariants. We differentiate with respect to $\bar z_s$,
$s=1,2$ and write the resulting differential equation in complex form
as
\begin{eqnarray}
  \label{eq:vf_compl_form}
  \dot z_1&=&\lambda z_1+c_1 z_1(|z_1|^2+|z_2|^2)+c_2z_1|z_2|^2
  +c_3i\bar z_1 z_2^2\\
  \dot z_2&=&\lambda z_2+c_2 z_2(|z_1|^2+|z_2|^2)+c_2z_2|z_1|^2
  +c_3i z_1^2\bar z_2.
\end{eqnarray}
For the last equivariant map we observe in real coordinates we take
$J\nabla I_{4,2}$. We get $\nabla I_{4,2}$ in the complex form by
computing the gradient with respect $\pt/\pt\bar z_j$ and multiplying
with $J$. The last operation is obtained by multiplying acting with
$[i,1]$ on the equation. This is the same as multiplying the $\bar
z$-gradient with $i$. Observe that in the case $G_3(m)$ we have to
take the real gradient of this invariant and obtain up to order $3$ a
fully gradient map. \\
Let us add one more remark: the number of cubic equivariant fields is
non increasing function in $m$. From the character formula we find,
that for $m=3$ we have three equivariant maps. So the vectorfield
given here is the cubic truncation in each case.\\
In a similar way we can construct the four sextic terms from
$I_2^3,I_{4,1}I_2,I_6$ and a Hamiltonian field from $I_{4,2}I_2$
(respectively non-Hamiltonian for $G_3(m)$. This listing is not
complete, the character formula predicts $9$ sextic equivariant
maps. Due to our method to require invariance for all $m$, we have
constructed the invariants and equivariants for the Lie groups
$G_j$. Therefore equilibria occur in circles, therefore they are not
hyperbolic. Of course this means that they could be destroyed by
adding higher order terms. Next we are looking at the fixed point
subspaces. Since these spaces are different for the three groups, we
discuss them one by one.
\begin{enumerate}
\item The case $G_1(m)$. \\
There is one isotropy type, a representative of this type is given by
the nontrivial element of order $2$: $\Z_2=\langle[j,je_4]\rangle$,
where $\langle\cdot\rangle $ denotes the group generated by the listed
elements. Its fixed point space can be explicitly computed by Lemma
\ref{lem:fixedpointspace_order_two}. We get 
\[\Meng{\alpha (1+e_4)+\beta (1+\bar
  e_4)j}{\alpha,\beta\in\R}=\Meng{\alpha e_8+\beta \bar e_8j}{\alpha,\beta\in\R}.\]
\item The case $G_2(m)$.\\
From the computation which we have given in the proof of Theorem
  \ref{series.thm:action} we can read off, two representatives for the
  two classes of isotropy subgroups. We have 
\[\Sigma_1=\langle[i,j]\rangle\subset F_2(m)\]
and
\[\Sigma_2=\langle [e_{2m}j,je_4]\rangle.\]
The fixed point spaces are given
\[\Fix{\Sigma_1}=\Meng{\alpha(1-ij)+\beta(i+j)}{\alpha,\beta\in\R}\]
and
\[\Fix{\Sigma_2}=\Meng{\alpha(1+\bar e_{2m}e_4)+\beta(e_{2m}+\bar
  e_4)j}{\alpha,\beta\in\R}.\]
\item The case $G_3(m)$.\\
We had seen that we have three conjugacy classes of groups of order
$2$, in each class we look at one representative, i.e. we look at
\[\Sigma_1=\langle[j,i]\rangle,\ \Sigma_2=\langle[j,j]\rangle,\
\Sigma_3=\langle[j,k]\rangle \]
with fixed point spaces given by Lemma \ref{lem:fixedpointspace_order_two} as
\begin{eqnarray*}
\Fix{\Sigma_1}&=&\Meng{\alpha(1-ji)+\beta(i+j)}{\alpha,\beta\in\R},\\
\Fix{\Sigma_2}&=&\Meng{\alpha+\beta j}{\alpha,\beta\in\R}\\
\Fix{\Sigma_3}&=&\Meng{\alpha(1-jk)+\beta(j+k)}{\alpha,\beta\in\R}\\
&=&\Meng{\alpha(1-i)+\beta(1+i)j}{\alpha,\beta\in\R}. 
\end{eqnarray*}
\end{enumerate}
The next observation comes from the equivariance with respect to the
Lie groups $G_j$. As a consequence the equations are equivariant with
respect to the action
\[\left(
  \begin{array}{c}
    z_1\\
    z_2
  \end{array}\right)\mapsto e^{i\phi} \left(
  \begin{array}{c}
    z_1\\
    z_2
  \end{array}\right).\]
The only fixed point of this group action is the origin and therefore
all equilibria are not isolated and therefore perturbations with
higher order terms might destroy all equilibria. Now we look at
restriction of the equation to the various fixed point spaces. We do
it again case by case.
\begin{enumerate}
\item $G_1(m)$\\
There is only one isotropy type with a nontrivial fixed point
space. The fixed point space is $\alpha e_8+\beta \bar e_8j$ and the
equation for equilibria reads (up to cubic order)
\begin{eqnarray*}
  0&=&\lambda \alpha e_8 +c_1\alpha e_8(\alpha^2+\beta^2)+c_2\alpha
  e_8\beta^2+c_3i\bar e_8\bar e_8^2\alpha\beta^2\\
  0&=&\lambda \beta \bar e_8+c_1\beta \bar e_8(\alpha^2+\beta^2)+c_2\beta
  \bar e_8\alpha^2+c_3i e_8^2 e_8\beta\alpha^2\\
\end{eqnarray*}
This gives the equation
\begin{eqnarray*}
  0&=&\alpha(\lambda  +c_1(\alpha^2+\beta^2)+(c_2 +c_3)
  \beta^2)e_8\\
  0&=&\beta(\lambda
  +c_1(\alpha^2+\beta^2)+(c_2-c_3)\beta\alpha^2)\bar e_8.
\end{eqnarray*}
We can look for solutions $(\lambda,\alpha,\beta)\in\R^3$ with
$\alpha\beta\not=0$. However then we get sign conditions on the
coefficients $c_1,c_2,c_3$. As a consequence we would not prove generic
bifurcation results. Let us look for solutions of the form
$(\lambda,0,\beta)$ with $\beta\not=0$. With the ansatz $\alpha=0$ the
first equation is identically satisfied, the second equation then leads
to
\[0=\lambda\beta+c_1\beta^3\]
which gives
\[\beta^2=-\frac{\lambda}{c_1}.\]
where we assume $0<|c_1|<\infty$ and we choose $\lambda$ such that the
term on the right hand side is positive. With the additional
assumptions $c_2+c_3\not=c_1$ we obtain the linearisation of the map
along the given branch as
\[\left(
  \begin{array}{cc}
    \lambda-\frac{c_2+c_3}{c_1}\lambda&0\\
    0& -2\lambda
  \end{array}\right),\]
which is regular and yields the persistence of this branch under
higher order perturbation.
 \item $G_2(m)$\\In this case we use the second fixed point space and
   we obtain precisely the same system as in the case before,
   therefore we obtain the same result.
 \item $G_3(m)$\\
Here we use the third fixed point space, again, we get the real
equation and therefore the same result
\end{enumerate}~
\endproof~
{\bf Proof of Theorem \ref{thm:persistance}.}
In each case we have constructed a branch in a generic third order
equation which is stable under higher order perturbation and hence the
result follows. 
\endproof
\section{GAP Computations}
\label{sec:gap}
\subsection{General remarks}

In this paper we have shown  that there are three infinite
series of groups of orders $48+32\mu$, $\mu\in\N$, which act
absolutely irreducibly 
on $\R^4$ and which have no odd-dimensional fixed point space. 
In this section we collect together some data obtained with
the computational group theory package GAP \cite{GAP},
using the groups in the Small Group Library.
For any of these groups, it is possible to
obtain its character table and hence determine whether it acts
irreducibly on $\R^4$. The subgroup lattice is obtained and then
character formulas are used to determine the isotropy subgroups and
the dimension of their fixed point spaces (see \cite{Matt-04} for further
details).  

This approach was applied to actions on $\R^4$
and also to actions on $\R^N$, for $4<N\leq 20$ with $N$ even. 
The results lead us to the following conjectures:
For dimensions $N=0\mod 4$, there are infinitely many groups acting
absolutely irreducibly on $\R^N$ that have no isotropy subgroups with
odd-dimensional fixed point spaces.
But for dimensions $N=2\mod 4$, there are no such groups. 

We have checked most of the groups of order up
to $1000$, however in the dimensions $4$, $8$ we did not look at the
groups of order $512$, due to the sheer number of
such groups: there are $10494213$ groups of order $512$. 
Even if this number of groups can be checked with a computer, there
are $49487365422$ groups of order $1024$ and this number of groups is
certainly out of reach for present day computers.


\subsection{The case $0\mod 4$}
The following tables gives the GAP numbers for finite groups of orders up
to 1000 (and some cases higher) which act
absolutely irreducibly and have only even dimensional fixed point
spaces. Some of these groups have several inequivalent representations
in these dimensions with the same properties. However we do not
provide this information.
The tables are based on computations on different computers
using the computer algebra package GAP.  

\subsubsection{$\R^4$}

The results of the GAP computations for actions on $\R^4$ are
summarised in Tables~\ref{tab:r4} and \ref{tab:r4a}, which list the
Small Group Library number of groups that act on $\R^4$ and have
no isotropy subgroup with an odd-dimensional fixed point space. 
Note that this list contains the
groups from the series $G_1(m)$, $G_2(m)$, $G_3(m)$ (compare
Table~\ref{tab:SmallGroupLibrary}), but also many other groups. 
Most of the groups in this table belong to the two-parameter families
\[
\pm \frac12[D_{2m}\times \overline{D}_{4n}],\qquad
\pm \frac14[D_{4m}\times\overline{D}_{4n}],\qquad
\pm [D_{2m}\times D_{2n}],
\]
in the notation of \cite{CS}. These families include $G_1(m)$,
$G_2(m)$, $G_3(m)$ in the case $n=2$. The remaining groups belong to
four of the one-parameter families in Table 4.1 of  \cite{CS},
\[
\pm [O \times D_{2n}] \mbox{ for } n=5,7 \mbox{ (480:969, 672:1053),} \]
\[
\pm \frac12[O \times D_{2n}] \mbox{ for } n=5,7,11,13,17  \mbox{ (240:106,
336:119, 528:91, 624:135, 816:102),} \]
\[
\pm \frac16[O \times D_{6n}] \mbox{ for } n=3,9,15  \mbox{ (144:32, 432:38, 720:106),} \]
\[
\pm [T \times D_{2n}] \mbox{ for } n=3,5,7,9,11,13,15,17  \]
\[\mbox{ (144:127,
240:108, 336:131, 432:262, 528:93, 624:147, 720:544, 816:104).}
\]
In most of these cases the isotropy subgroups are isomorphic to
$\ZZ{2}$. But in some cases, for example the group 144:127, there is
an isotropy subgroup isomorphic to $\ZZ{3}$.

\begin{table}
\begin{tabular}{|c|l|}
\hline 
Order&Groups\\
\hline 
$48$&$15$ $17$ $41$\\
\hline
$80$&$15$ $17$ $42$\\
\hline
$96$&$35$ $33$ $126$\\
\hline
$112$&$14$ $16$ $34$\\
\hline
$120$&$10$ $12$ $13$\\
\hline
$144$&$16$ $18$ $32$ $44$ $127$\\
\hline
$160$&$35$ $33$ $140$\\
\hline
$168$&$14$ $16$ $17$\\
\hline
$176$&$14$ $16$ $34$\\
\hline
$192$&$78$ $80$ $478$\\
\hline
$208$&$15$ $17$ $42$\\
\hline 
$224$&$32$ $34$ $114$\\
\hline
$240$&$14$ $15$ $19$ $21$ $76$ $78$ $106$
      $108$\\&  $130$ $134$ $182$\\
\hline
$264$&$7$  $9$  $10$\\
\hline
$272$&$15$ $17$ $43$\\
\hline
$280$&$9$  $11$ $12$\\
\hline
$288$&$33$ $35$ $129$\\
\hline
$304$&$14$ $16$ $34$\\
\hline
$312$&$17$ $19$ $20$\\
\hline
$320$&$77$ $79$ $546$\\
\hline
$336$&$30$ $31$ $35$ $37$ $101$ $103$ $119$ $131$
      \\& $142$ $146$ $201$\\ 
\hline
\end{tabular} 
\begin{tabular}{|c|l|}
\hline
Order&Groups\\
\hline
$352$&$32$ $34$ $114$\\
\hline
$360$&$9$ $10$ $13$\\
\hline
$368$&$14$ $16$ $34$\\
\hline
$384$&$183$ $185$ $1959$\\
\hline
$400$&$16$ $18$ $42$\\
\hline
$408$&$9$ $11$ $12$\\
\hline
$416$&$33$ $35$ $140$\\
\hline
$432$&$16$ $18$ $38$ $50$ $262$\\
\hline
$440$&$19$ $21$ $22$\\
\hline
$448$&$76$ $78$ $453$\\
\hline
$456$&$14$ $16$ $17$\\
\hline
$464$&$15$ $17$ $42$\\
\hline
$480$&$14$ $15$ $19$ $21$ $186$ $188$ $347$ $351$\\& $884$ 
    $969$\\
\hline
$496$&$14$ $16$ $34$\\
\hline
$504$&$14$ $15$ $18$\\
\hline
$520$&$13$ $15$  $16$\\ 
\hline
$528$&$12$ $13$  $17$ $19$ $74$ $76$ $91$\\& $93$ $101$ $105$ $153$\\ 
\hline
$544$&$33$ $35$ $141$\\ 
\hline
$552$&$7$ $9$ $10$\\ 
\hline
$560$&$14$ $15$ $19$ $21$ $75$ $77$ $114$\\& $118$ $163$\\
 \hline
\end{tabular}
\caption{Small Group Library numbers of groups acting on $\R^4$ that
  have no odd-dimensional isotropy subgroups.  The groups of order
  $512$ have not been checked.\label{tab:r4}} 
\end{table}

\begin{table}
\begin{tabular}{|c|l|}
\hline 
Order&Groups\\
\hline 
$576$&$78$ $80$ $481$\\
\hline
$592$&$15$ $17$ $42$\\
\hline
$600$&$9$  $10$ $13$\\
\hline
$608$&$32$ $34$ $114$\\
\hline
$616$&$7$ $9$ $10$\\
\hline
$624$&$33$ $34$ $38$ $40$ $104$ $106$\\&  $135$ $147$  $172$ $176$
$231$\\ 
\hline
$640$&$180$ $182$ $2258$\\
\hline
$656$&$15$ $17$ $43$\\
\hline
$672$&$55$ $56$ $60$ $62$ $247$ $249$ $415$\\& $419$ $984$ $1053$\\
\hline
$680$&$13$ $15$ $16$\\
\hline
$688$&$14$ $16$ $34$\\
\hline
$696$&$9$ $11$ $12$\\
\hline
$704$&$76$ $78$ $451$\\
\hline
$720$&$12$ $13$ $18$ $22$ $77$ $79$ $106$\\& $127$ $134$ $182$ $544$\\
\hline
$728$&$9$ $11$ $12$\\
\hline
$736$&$32$ $34$ $114$\\
\hline
$744$&$14$ $16$ $17$\\
\hline
$752$&$14$ $16$ $34$\\
\hline
$760$&$9$  $11$  $12$\\
\hline
$784$&$15$ $17$ $34$\\
\hline
\end{tabular}
\begin{tabular}{|c|l|}
\hline
Order&Groups\\
\hline
$792$&$7$ $8$ $11$\\
\hline
$800$&$33$ $35$ $140$\\
\hline
$816$&$14$ $15$ $19$ $21$  $76$ $78$ $102$\\& $104$  $126$  $130$ $178$\\
\hline 
$832$&$77$ $79$ $546$\\
\hline 
$840$&$58$ $60$ $61$ $72$ $74$ $75$ $79$ $81$ $82$\\
\hline
$848$&$15$ $17$ $42$ \\ 
\hline
$864$&$34$ $36$  $135$\\
\hline
$880$&$40$ $41$ $45$ $47$ $101$ $103$  $154$ $158$ $203$\\
\hline
$888$&$17$ $19$ $20$\\
\hline
$896$&$179$    $181$  $1903$  \\
\hline
$912$&$30$  $31$  $35$ $37$  $101$ $103$  $118$ $130$  $141$\\& $145$  $200$ \\
\hline
$920$&$9$   $11$  $12$   \\
\hline
$928$&
$33$ $35$  $140$ \\
\hline
$936$&$17$ $18$  $21$   \\
\hline
$944$&$14$ $16$ $34$   \\
\hline
$952$& $9$   $11$  $12$    \\
\hline
$960$&$14$  $15$  $19$  $21$  $535$  $537$  $1029$  $1033$\\&  $5215$  \\
\hline
$976$&$15$   $17$  $42$    \\
\hline
$984$&$9$    $11$  $12$   \\
\hline
$992$&
$32$    $34$  $114$  \\
\hline
\end{tabular}\\[1cm]
\caption{Small Group Library numbers of groups acting on $\R^4$ that
  have no odd-dimensional isotropy subgroups (continued).\label{tab:r4a}} 
\end{table}
All the groups which we have discussed here are subgroups of
$\SO{4}$. It is easy to see, that groups with elements in
$\OO{4}\setminus \SO{4}$ do have odd dimensional fixed point
spaces. If $g$ is such an element, its determinant is $-1$, so $-1$ is
an eigenvalue of multiplicity $1$ or $3$. In the second case, $1$ is
an eigenvalue of multiplicity one. In the first case, $1$ is an
eigenvalue of multiplicity $1$ or $3$. So in any case the group
generated by $g$ has an odd dimensional fixed point space. 

\subsubsection{$\R^8$}
Table ~\ref{tab:r8} is a similar list for groups acting absolutely
irreducibly on $\R^8$.  The case of groups of order $512$ has not been
checked. In the case of 
groups of order $768$ we have the complete answer, the list presented
here is only part of the list we have obtained so far. 

\begin{table} 
\begin{tabular}{|c|l|}
\hline
Order&Groups\\
\hline
$160$&$82$  $85$  $208$\\
\hline
$192$&$36$  $308$  $310$  $312$  $313$  $758$  $761$  $762$  $804$  $990$  $1333$  $1336$  $1337$  $1394$  $1396$  $1484$  $1527$  \\
\hline
$240$&$96$  $99$  $101$  \\
\hline
$288$&$382$   $383$  $433$   $572$  $573$  $582$  $583$
$586$  $587$  $589$  $593$  $596$  $597$  $598$   $937$  $964$
$966$  $968$  \\
\hline
$320$&
$35$  $242$ $244$  $245$ $247$ $266$
$376$ $378$ $380$ $381$ $826$  $829$ $830$ $872$ $1078$ 
\\& $1079$ $1122$ $1123$ $1446$  $1449$ $1450$  $1507$ $1509$  $1598$ $1625$\\
\hline
$384$&$114$  $117$   $126$  $128$  $346$  $349$ $1744$ $1748$  $1834$  $1836$ $1842$ 
$1844$  $1847$  $1848$  $3570$ 
\\&  $3668$ $3670$ $3673$  $3674$   $4093$  $4095$ $4099$  $4698$  $4700$  $5713$ 
$5715$   $12858$  $12874$  \\&  $12878$
$13519$ $13522$  $14612$   $14615$
$14616$   $16583$  $16592$   $16596$
$19786$  \\
\hline
$416$&$82$ $85$ $207$ \\
\hline
$448$&$34$ $283$ $285$
$287$ $288$ $733$ $736$ 
 $737$ $779$  $1227$ $1230$ $1231$  $1288$ $1290$ $1382$\\ 
\hline
$480$&$227$ $228$  $233$  $234$  $249$  $250$  $553$ -- $557$  $563$   $567$  $568$  $571$
$572$  $574$  $576$ --  $580$\\&  $582$  $588$  $591$
$592$ $595$ -- $597$  $599$ 
$959$  $961$  $964$  $969$   $970$   $973$\\&  $990$
 $992$  $1007$  $1100$  $1102$  
$1103$ $1105$  $1108$  $1109$  $1111$ \\
\hline
$544$&
$83$   $86$  $216$ \\
\hline
$560$&
$89$   $92$  $94$  \\
\hline
$576$&
$37$   $311$  $313$ $315$  $316$ $762$  $765$  $766$  $808$ $1065$  $1415$   $1758$  $1761$  $1762$ $1819$  \\&  $1821$
$1906$  $1907$  $1927$  $1938$  $1940$  $1946$  $1947$  $1949$  $2078$ $2080$ $2099$  $2905$    \\& $2906$  $2915$  $2916$ $2919$  $2920$   $2922$  $2926$  $2929$  $2930$  $2931$  $3390$   $3522$  $3525$\\&  $3528$  $3540$  $4982$  $4983$ 
$5025$  $5106$  $5107$  $5117$  $5216$   $5228$  $5263$  $5297$  $5547$   \\&  $5599$  $5601$ $5669$ $5670$  $5713$  $5714$  $6643$  $6644$  $6647$  $6649$  $6653$  $6979$  $7201$  \\& $7205$  $8273$  $8317$
 $8330$  $8332$  $8338$  $8470$ $8507$ $8526$  $8571$  \\
\hline
$624$&
$125$   $128$  $130$  \\
\hline
$640$&
$111$ $114$  $123$ 
$125$  $343$  $346$   $760$
$763$  $765$  $766$   $829$  $832$ $835$  $839$ $2043$  $2047$ 
 \\& $2133$  $2135$ $2141$
$2143$   $2146$  $2147$
$3869$  $3967$ 
  $3969$  $3972$  $3973$  $4392$  $4394$  \\&  $4398$  $4997$  $4999$ $6218$  $6219$  $6219$ 
$6722$ $6912$
$6916$ $6929$
$14095$  $14111$
  \\&$14115$
$15852$   $15853$ $14756$
$14759$ $15849$
$17820$   $17829$
$17833$  $19519$
$21193$   \\ 
\hline
$672$&
$621$ -- $625$  $631$  $635$ $636$ $640$
$642$    $644$  $645$ -- $648$  $650$   $656$  $659$ $660$
$663$  \\ & $664$ --  $667$  $1053$  $1054$  $1057$ $1153$  $1155$ $1156$  $1158$   $1161$  $1162$ $1164$  \\
\hline
$704$&
$34$   $281$  $283$  $285$  $286$  $731$  $734$  $735$  $777$  $1224$
$1227$  $1228$  $1285$  $1287$  $1376$  \\
\hline
$720$&
$96$   $98$  $101$   $450$  $452$   $457$  $459$   $460$  $462$  $475$  $476$  $481$  $482$ $490$  $495$   $496$    \\&$501$   $502$  $504$  $509$ $510$ $513$  $517$  $520$  $521$   $523$  $524$ \\
\hline
$768$& $57401$ $57403$ $57429$ $80778$ $82966$ $82967$ $82970$ $83806$
$83807$ $83820$ $83821$\\ &$89833$ $89899$ $90043$ $90249$  $90250$ 
$90252$ $90255$ $90259$  $90262$ $90263$ $90264$ \\
\hline
\end{tabular}
\caption{As Table~\ref{tab:r4} but for actions on $\R^8$. \label{tab:r8}}
\end{table} 

\subsubsection{$\R^{12}$}

Table~\ref{tab:r12} is a list of groups acting on $\R^{12}$, with the
same property. 

\begin{table}
\begin{tabular}{|c|l|}
\hline
Order&Groups\\
\hline
$336$&$18$ $20$ $128$ $134$\\
\hline
$432$&$83$ $85$ $153$ $155$ $161$ $163$ $245$ $248$ $261$ $269$
      $295$ $369$ $371$\\ 
\hline
$504$&$69$ $71$ $72$\\
\hline
$576$&$184$ $1399$ $1401$ $1990$ $2010$ $4976$ $4988$ $5060$ $5061$
      $5065$ $5067$ $5531$ \\& $5533$ $5584$ $5586$ $8307$ $8310$
      $8464$ $8483$ $8499$ $8510$\\
\hline
$624$&$19$    $21$  $144$   $150$  \\
\hline
$672$&
$36$    $38$  $337$  $337$  \\
\hline
$720$&$415$  $418$  \\
\hline
$840$&
$18$   $20$    $21$  \\
\hline
$864$&
$220$   $222$  $428$   $430$  $436$    $438$  $662$   $822$  $1157$  $1195$  $1206$  $1532$   \\& $1535$   $2195$  $2222$  $2271$ $2487$  \\
\hline
$912$&
$18$    $20$  $127$   $133$  \\
\hline
$936$&
$76$  $78$   $79$  \\
\hline
$960$&
$789$    $809$  $5713$  $5714$  $5769$    $5770$  $5774$   $5776$  $6329$   $6331$  $6382$   $6384$ \\&   $10946$ $10949$   $11098$  $11117$ $11133$  $11144$  \\
\hline
$1008$&
$224$   $225$  $229$ $231$  $286$   $288$  $522$  $524$  $532$    $536$  $664$   $682$  $881$  \\
\hline
$1080$&
$92$    $93$  $96$  $99$  $100$   $103$  $145$   $147$  $149$   $282$  \\
\hline
$1152$& 
$153939$ $153941$ $153959$ $153960$  $153963$ $153969$
     $154100$ $154102$ $154147$ \\& $154149$
     $154375$  $154498$ $154503$  $154506$ $154507$ $154590$  $154592$  $154596$ \\& $154598$ $154690$ $154960$
     $154961$ $154964$ $154970$
     $155100$ $155102$  $155172$ \\& $155191$ $155194$ $155356$ $155358$
      $155809$ $155823$ $156072$
     $156207$  $156208$  \\&$156214$ $156214$ $157025$
    $157340$ $157585$ $157644$ \\
\hline
\end{tabular}
\caption{As Table~\ref{tab:r4} but for actions on $\R^{12}$. \label{tab:r12}}
\end{table} 

\subsubsection{$\R^{16}$}
A similar list for groups acting on $\R^{16}$ is given in
Table~\ref{tab:r16}. 
The gap numbers of all groups of order $768$ have been determined,
but the number is too large to present all of them here. 

\begin{table}
\begin{tabular}{|c|l|}
\hline
Order&Groups\\
\hline
$576$&$5153$  $8369$  \\
\hline
$640$&$653$ $657$ $660$ $915$ $917$ $6009$ $6012$ $6014$ $6016$ $6948$
$6950$ $6951$ $7102$ $19529$\\& $19534$ $19535$ $19641$ $19644$ $21497$\\ 
\hline
$768$&$57413$ $57415$  $79717$  $79730$  $80103$  $80107$  $80108$ $80111$  $80408$  $80412$  $8055
7$  \\
&\dots\\
&$1045835$ $1045841$  $1045859$  $1045863$ $1045866$\\
\hline
$960$&$6130$  $6131$  $6133$  $6134$  $6135$  $6136$  $6137$  $6144$  $6148$  $\
6149$  $6152$  $6153$  $6156$\\&  
$6158$  $6159$  $6160$  $6161$  $6166$  $6167$ $6169$  $6172$  $6173$
$6176$  $6177$  $10895$  \\&  
$10897$  $10900$ $11035$  $11037$  $11040$  $11043$  $11046$  $11048$  $11049$  \\
\hline
\end{tabular}
\caption{As Table~\ref{tab:r4} but for actions on $\R^{16}$. \label{tab:r16}}
\end{table} 

\subsubsection{$\R^{20}$}
Here we present the list for the groups acting on $\R^{20}$.
\begin{table}\begin{center}
\begin{tabular}{|c|l|}
\hline
Order&Groups\\
\hline
$880$& $18$ $20$ $121$\\
\hline
$1320$&$17$ $19$ $20$\\
\hline 
$1760$&$36$  $38$  $329$  \\
\hline
\end{tabular}
\caption{As Table~\ref{tab:r4} but for actions on $\R^{20}$. \label{tab:r20}}
\end{center}
\end{table}


\subsection{The case $2\mod 4$}
We have already remarked that in the case of $\R^2$ there is no
absolutely irreducible representation without an odd dimensional fixed
point space.
We have checked all groups up to the following orders in the various
dimensions and have not found any groups which act absolutely
irreducibly and have no odd dimensional fixed point spaces.
\[\begin{array}{|c|c|}
\hline
{\rm dimension}&{\rm Order} \\
\hline
2&\infty\\
6& 1013\\
10&999 \\
14&1007\\
18&1151\\
\hline 
\end{array}
\] 
At this point we mention a recent result by Ruan \cite{Ruan10}: in
dimension $6$ all solvable groups which act absolutely irreducibly have
an odd dimensional fixed point space.  
\section{Characters and Invariant Theory}\label{characters}
In Sattinger \cite{Sat1} we find a formula which allows to compute the
vector space dimensions of the space of invariant polynomials for a group
action of a given degree.
Consider the action of a compact Lie group $G$ on some finite dimensional
space $V$. We write $C^\infty_G(V)$ for the $G$ invariant smooth functions. It
is well known that they form an algebra which is finitely generated by invariant
polynomials. Therefore we are interested in the space of homogeneous invariant
polynomials of some given degree $d$. The dimension of this space is denoted by
$c_d$ and in a similar way we write $C_d$ for the dimension of the space of
homogeneous equivariant polynomial maps for $V\to V$. Sattinger \cite{Sat1}
defines the quantities  
\[\chi_{(d)}(g)=\displaystyle\sum_{\displaystyle\sum_{j=1}^d j\,i_j=d}
\displaystyle\frac{\chi^{i_1}(g)\chi^{i_2}(g^2)\chi^{i_3}(g^3)\cdots\chi^{i_d}(g^d)}{1^{i_1}\,i_1!\,2^{i_2}\,i_2!\,\cdots\,d^{i_d}\,i_d!}\] 
and obtains the following representations for $c_d, C_d$:
\begin{equation}
  \label{eq:inv}
c_d=\intl{G}{}{\chi_{(d)}(g)\;dg}  
\end{equation}
and
\begin{equation}
  \label{eq:equ}
  C_d=\intl{G}{}{\chi_{(d)}(g)\chi(g)\;ds}.
\end{equation}
We obtain for (following Sattinger \cite{Sat1})
\[\chi_{(2)}=\frac12\left(\chi(g^2)+\chi^2(g)\right).\]
For the next values we derive the following expressions (using that
\[i_1+2i_2+3i_3=3\]
leads to the choices $(3,0,0), (1,1,0), (0,0,1)$ for $(i_1,i_2,i_3)$
and
\[i_1+2i_2+3i_3+4i_4=4\]
to $(4,0,0,0), (1,0,1,0), (2,1,0,0),(0,2,0,0), (0,0,0,1)$ for $(i_1,i_2,i_3,i_4)$)
\[\chi_{(3)}=\frac1{3!}\chi^3(g)+\frac12\chi(g)\chi(g^2)+\frac{1}{3}\chi(g^3).\]
and
\[\chi_{(4)}=\frac1{4!}\chi^4(g)+\frac13\chi(g)\chi(g^3)+\frac14\chi^2(g)\chi(g^2)+\frac1{2^2\,2!}\chi^2(g^2)+\frac1{4}\chi(g^4).\]

\addcontentsline{toc}{chapter}{Literaturverzeichnis}
\bibliography{all,bif_a,bif_b,bif_c,bif_d,bif_f,bif_g,bif_h,bif_k,bif_m,bif_v,bifurcat,pde,mech,a,vorlesung,fa} 
\bibliographystyle{ab1a}


\end{document}
